\documentstyle [11pt]{article}
\newfont{\bbb} {msbm10}
\newcommand{\Bbb}[1]{\mbox{\bbb#1}}

\newcommand{\h}{\breve{H}}
\newcommand{\bX}{\bar{X}}

\newcommand{\cU}{{\cal{U}}}
\newcommand{\cN}{{\cal{N}}}
\newcommand{\sm}{\setminus}
\newcommand{\sbs}{\subset}
\newcommand{\ra}{\rightarrow}
\newcommand{\p}{\partial}
\begin{document}
\title{ Cocompact Proper CAT(0) Spaces}
\author{Ross Geoghegan and Pedro Ontaneda\thanks{
The second-named author was supported in part 
by a research  grant from CAPES, Brazil.}}
\maketitle

\begin{abstract}
This paper is about geometric and topological properties of a proper 
CAT(0) space $X$ which is cocompact - i.e. which has a compact 
generating domain with respect to the full isometry group.  It is 
shown that geodesic segments  in $X$ can ``almost" be extended to 
geodesic rays.  A basic ingredient of the proof of  this geometric 
statement is the topological theorem  that there is a top dimension 
$d$ in which the compactly supported integral cohomology of $X$ is 
non-zero.  It is also proved that the boundary-at-infinity of $X$ 
(with the cone topology) has Lebesgue covering dimension $d-1$.  It 
is not assumed that there is any cocompact discrete subgroup of the 
isometry group of $X$; however,  a corollary for that case is that 
``the dimension of the boundary" is a quasi- isometry invariant of 
CAT(0) groups.  (By contrast, it is known that the topological type 
of the boundary is not unique for a CAT(0) group.)
\end{abstract} 

\section{Statement of Theorems}

A CAT(0) space is a geodesic metric space $(X, d_X)$  whose geodesic
triangles are ``no fatter than" the corresponding comparison 
triangles in
the Euclidean plane.  A general reference for facts
about CAT(0) spaces used here is \cite{BH}.  We will usually suppress
$d_X$ referring just to $X$.   Such a space $X$ is {\it proper}
if all closed balls are compact, and is {\it cocompact} if there is
a compact generating domain for the full isometry group of $X$,
i.e. there is a compact set $C \subset X$ such that the sets $\{h(C)
\mid h$ is an isometry of $X\}$ cover $X$.   In particular, a proper 
CAT(0) space $X$
has a compact {\it boundary}, $\partial _{\infty} X$, namely  the set
of asymptoty classes of geodesic rays in $X$ with the ``cone 
topology".
Equivalently, picking a base point $p \in X$ one can regard $\partial
_{\infty} X$ as the set of geodesic rays starting at $p$ endowed with
the compact-open topology.\\

The boundary of a proper CAT(0) space can be infinite-dimensional 
(see Example 4 below) but there is the following theorem of E. L. 
Swenson:\\

\noindent {\bf {Theorem 0.}} (\cite{Sw}; Theorem 12). {\it If the 
proper CAT(0) space $X$ is
cocompact then $\partial _{\infty} X$ has finite (Lebsegue covering) 
dimension.}\\

In this paper we sharpen Swenson's theorem by showing that 
cocompactness implies more.  As usual $H_c^*(X)$ denotes the integral 
cohomology of $X$ with compact supports.  For cocompact $X$ we prove 
{\it (i)} that there is a highest dimension $d$ for which $H_c^d(X)$ 
is non-zero, and {\it (ii)} that $\partial _{\infty} X$ has dimension 
$d-1$.   While this is of interest in its own right, our initial 
motivation came from a wish to prove that every cocompact non-compact 
proper CAT(0) space is almost geodesically complete (Theorem 5, 
below).  By previous work of the second-named author \cite {O} this 
follows from{\it (i)}.  We know of no geometric proof of this 
geometric statement, i.e. a proof of almost geodesic completeness 
for proper cocompact CAT(0) spaces, which does not use algebraic 
topology.\\

Before stating our results in detail we need some definitions.\\

Two metric spaces {\it have the same bounded homotopy type} if there 
are maps (i.e. continuous functions) from each to the other 
such that either composition is homotopic to the appropriate identity 
map by a homotopy which only moves points by a bounded amount, i.e. 
there is a number $s\geq 0$ such that the image of ${point} \times 
[0,1]$ under either of the homotopies has diameter at most $s$.  
When that happens the maps in question are called {\it bounded 
homotopy equivalences}.  If the metric spaces are proper (i.e all 
closed balls are 
compact) then these must be proper; moreover,  ``having  the same 
bounded homotopy type" implies ``having the same proper homotopy 
type" and a bounded homotopy equivalence is a 
proper homotopy equivalence.  In particular,  $H^*_c$ is a bounded 
homotopy invariant.\\

When $K$ is a countable locally finite simplicial complex,  $|K|$ 
will be understood to carry the {\it unit metric} unless we say 
otherwise: by this we mean that each simplex is isometric to a 
standard Euclidean simplex whose edges have length 1, and the 
distance between points of $|K|$ is the greatest lower bound of the 
lengths of all piecewise linear paths, each piece in a simplex, 
joining them.  Thus $|K|$ is a proper geodesic metric space.\\

\noindent {\bf Theorem 1.} {\it Let $X$ be a cocompact proper CAT(0) 
space.  (i) There is a finite-dimensional  countable locally finite 
simplicial complex $K$ such that $X$ and $|K|$ have the same bounded 
homotopy type.    (ii)  There is a number $d$ such that $H_c^d(X)$ is 
non-trivial while $H_c^i(X)$ is trivial for all $i > d$.}\\

A version of Theorem 1 appears in \cite{O} but with the additional 
hypothesis that for some cocompact group of isometries $\Gamma$ the 
orbits of the $\Gamma$-action are discrete; see  footnote on page 209 
of \cite{BH}.  One of the points of the present paper is to remove 
such discreteness hypotheses.\\

It is well known (see for example \cite{G_1}, \cite{G_2} or 
\cite{GM}) that $H_c^n(X)$ is isomorphic to the reduced integral 
$(n-1)$-dimensional {\v C}ech cohomology of $\partial _{\infty} X$.   
Thus we have:\\

\noindent {\bf Corollary 2.} {\it $\partial _{\infty} X$ has 
non-trivial reduced integral {\v C}ech cohomology in dimension $d-1$ 
and trivial 
integral {\v C}ech cohomology in all higher dimensions.}\\

\noindent {\bf  Example 3.}   Let $Q$ denote the Hilbert Cube in
Hilbert space, i.e. the set of sequences of real numbers whose $n$th
entry lies in the closed interval $[2^{-n}, 2^n]$, with
metric coming from the Hilbert space norm.  Hilbert space is obviously
CAT(0), so $Q$, being a convex subset, is also CAT(0). The boundary
of $\Bbb{R}$ is homeomorphic to $S^0$, but the boundary of $\Bbb{R}
\times Q$ is also homeomorphic to $S^0$, a simple illustration of how
infinite-dimensional CAT(0) spaces can have finite-dimensional 
boundaries.
Similarly $\partial {{\Bbb{R}}^{n} \times Q}$ is an $n- 1$-sphere. 
These
are cocompact examples with $d=n$.\\

\noindent {\bf Example 4.} The proper CAT(0) spaces ${\Bbb{R}}_{\geq 
0} \times {\Bbb{R}}^{n}$ and 
``the open cone on Q'' (made into a CAT(0) space in the
way described in  \cite{BH} II.3.14)  have 
contractible boundaries homeomorphic to $B^n$ and $Q$ respectively so 
Corollary 2 implies that they  are not cocompact.  (Indeed, by 
Theorem 1(ii), the
underlying topological spaces do not admit proper cocompact CAT(0) 
metrics because in this case the cohomology with compact supports is 
trivial.)\\

Our main application of Theorem 1 is to geometry.   A CAT(0) metric 
space
$X$ is {\it almost geodesically complete} if there is a number $r \geq
0$ such that for any points $a$ and $b$ of $X$ there is a geodesic ray
$\gamma : [0, \infty ) \ra X$ with $\gamma (a) = 0$ whose image meets
the open ball about $b$ of radius $r$.  The term ``almost extendible"
is also used. Obviously, the spaces in Example 3  have this property
while those in Example 4 do not.  An example which has this property 
even
though not every geodesic segment can be extended to a geodesic ray 
is the
graph consisting of  $\Bbb{R}$ together with, for each $n \in 
\Bbb{Z}$,
a copy of the closed unit interval glued at its $0$-point to $n\in  
\Bbb{R}$.\\

Theorem A of \cite{O} says that {\it if there were a non-compact 
cocompact proper CAT(0) space $X$ which is not almost geodesically 
complete then $H_c^*(X)$ would have to be trivial}.  Combining this 
with Theorem 1 we get:\\

\noindent {\bf Theorem 5.} {\it Every non-compact cocompact proper 
CAT(0) space is almost geodesically complete.}\\

As explained above, this was proved in \cite{O} under the additional 
hypothesis that there exists a  cocompact group of isometries with 
discrete orbits.\\

We remark that Theorem H of \cite{BG} requires almost geodesic 
completeness as a hypothesis on the non-compact cocompact proper 
CAT(0) space under consideration.  Theorem 5 shows that this 
hypothesis is redundant.\\

Corollary 2 gives $d-1$ as the sharp upper bound for non-vanishing  
\v {C}ech cohomology, but we can prove a better statement:\\

\noindent {\bf Theorem 6.} {\it Let $X$ be a cocompact proper CAT(0) 
space.  Then the (Lebesgue covering) dimension of $\partial _{\infty} 
X$ is $d-1$, where $d$ is as in Theorem 1.}\\

We remark that Theorem C of \cite{K} might be considered to be an 
analog
of our Theorem 6 for $\partial _{\infty} X$ when it is equipped with 
the Tits metric topology.\\

A {\it CAT(0) group} is a group $\Gamma$ 
which can act {\it geometrically} (i.e. properly discontinuously and 
cocompactly by isometries) on some proper CAT(0) space $X$.  Then 
$H_c^*(X)$ is isomorphic to $H^*(\Gamma, \Bbb{Z}\Gamma)$ (see 
Exercise VIII.7.4 of  \cite{Br}), so the number $d$ in Theorem 1 
depends only on $\Gamma$, hence Theorem 6 implies that the dimension 
of $\partial _{\infty} X$ depends only on $\Gamma$.  This is of 
interest because Croke and Kleiner \cite{CK} have given examples to 
show that the homeomorphism type of $\partial _{\infty} X$ is not an 
invariant of $\Gamma$.  Indeed, if $\Gamma_1$ and $\Gamma_2$ are 
quasi-isometric CAT(0) groups then $H^*(\Gamma_1, \Bbb{Z}\Gamma_1)$ 
and  $H^*(\Gamma_2, \Bbb{Z}\Gamma_2)$ are isomorphic \cite{Ger}, so 
we have:\\

\noindent {\bf Corollary 7.} {\it The Lebesgue covering dimension of 
the boundary of proper CAT(0) spaces on which CAT(0) groups act 
geometrically is a quasi-isometry invariant (of CAT(0) groups)}.\\

Of course, if $\Gamma$ acts cocompactly as covering transformations
on $X$ then $\Gamma$ has finite cohomological dimension so the number
$d$ of Theorem 1 exists for group theoretic reasons, and in that case
Corollary 7 follows from the proof of Corollary 1.4 of \cite{BM}, as 
is
pointed out in \cite{Be}.

\section{Proofs.}

We are to prove Theorems 1 and 6.\\

\noindent {\bf Proof of Theorem 1(i).} Let $E\sbs X$ be maximal with
respect to the property that if $x,y\in E$ and  $x\neq y$, then $d_X
(x,y)\geq 1$.  The family $\cU =\{ B_{X}(x,1) \mid  x\in E\}$ is an 
open
cover of $X$, where $B_{X}(x,1)$ denotes the open ball of radius 1.
Let $K$ be the nerve of this cover. Then $K$ is clearly countable and
it is locally finite because the cover $\cU$ is star finite (i.e. each
member of $\cU$ meets only finitely many others).\\

Suppose $K$ is not finite-dimensional.  Then for all natural numbers 
$m$ we would have dim$K$ $\geq m$.  Thus for each $m$ there would be 
points $x_0^m, ..., x_m^m \in E$ such that $1\leq d_X(x_i^m, x_j^m)< 
2$ when $i \neq j$.  Let $C$ be a compact generating domain for $X$.  
For each $m$ there is an isometry $h_m$ such that $h_m(x_0^m) \in 
C$.  Then for each $i$ the point  $h_m(x_i^m)$ lies in $C_1$, the 
1-neighborhood 
of $C$.  By induction, we can then pick subsequences of $\Bbb{N}$, 
each a subsequence of its predecessor, so that for each $n$ the 
sequence 
$(h_{n+k}(x_n^{n+k}))$ converges to some point $y_n \in C_1$.  
The resulting sequence $(y_n)$ consists of points which are pairwise 
at least 1 apart while all lying in the compact set $C_1$, a 
contradiction.\\

Because $X$ is CAT(0) and $\cU$ consists of open balls, it follows
that all non-empty intersections of members of $\cU$ are contractible.
It is well known that this implies that $X$ is homotopy equivalent
to $|K|$ (the weak topology coincides with the unit metric topology
since $K$ is locally finite), and indeed the proof of this fact in \S
5 of \cite{W}  shows that, with the metric we have chosen, the mutually
homotopy inverse maps in both directions given in that proof  are bounded
homotopy equivalences which are bounded homotopy inverse to one another.
(Indeed, this also follows from the proof of Lemma 7A.15 on page 129
of \cite{BH}).  q.e.d.\\

When $Z$ is a metric space $B_Z(p,r)$ denotes the open ball of radius 
$r$ and center $p\in Z$.   We say that $Z$ is {\it uniformly 
contractible} if for every $r>0$ there is $s>0$ such that for every 
$p \in Z$, $B_Z(p,r)$ contracts in $B_Z(p,s)$.  We say  $H^i_c(Z)$ 
{\it is uniformly trivial} if $Z$ has the following property:\\

\noindent {\it For every $r>0$ there is $s>0$ such that whenever a 
$i$-cocycle $z$ has compact support contained in a ball $B_Z(p,r)$, 
then $z$ cobounds a cochain whose compact support  lies in 
$B_Z(p,s)$.}  \\

Recall that when $U \subset V \subset Z$ with $U$ and $V$ both open 
in $Z$ then the inclusion map $\iota: U \ra V$ induces a 
homomorphism $\iota_*: H_c^*(U)\ra H_c^*(V)$; see Remark 26.2 of 
\cite{GH}.\\ 

\noindent {\bf Proposition 2.1.} {\it Let  $X$ be a proper cocompact 
CAT(0) space. If $H_c^i(X)$ is trivial, then $H_c^i(X)$ is uniformly 
trivial.}\\

\noindent {\bf Proof.}
Let $r>0$ and let $x\in X$.\\

\noindent {\bf Claim.} {\it There is $r'\geq r$ (depending only on 
$r$) such that
$\iota_*(H_c^i(B_X(x,r)))$ is a finitely generated subgroup of 
$H_c^i(B_X(x,r'))$,
where $\iota :B_X(x,r)\ra B_X(x,r')$ is the inclusion.}\\

To prove the Claim, let $K$ be as in the proof of Theorem 1(i) and let
$f:X\ra |K|$, $g: |K|\ra X$ be such that $gf$ is bounded homotopic
to the identity map.  Then there is a compact set $C$ containing 
$B_X(x,r)$ such that, for any cocycle $z$ with compact support lying 
in $B_{X}(x,r)$, $z$ and $f^*g^*z$ are compactly cohomologous in $C$. 
\\

Let $L$ be a finite subcomplex of $K$ such that:\\

${\bf(*)}\hspace{1.2in} g^{-1}(C)\subset int\ |L|$.\\

Choose $r'$ such that $g(|L|)\subset B_X(x,r')$.   It follows that 
$C\subset g(int\, |L|) \subset B_X(x,r')$.  Hence, if $z$ is a cocycle
with compact support contained in $B_X(x,r)$, then $z$ is compactly
cohomologous to $(gf)^*z=f^*g^*z$ in $B_X(x,r')$.  
By (*)
$v:=g^*z$ is a cocycle with compact support contained in $int\, |L|$.
Consequently every cocycle with compact support contained in 
$B_X(x,r)$
is compactly cohomologous in $B_X(x,r')$ to a cocycle of the form 
$f^*v$
where $v$ is a cocycle with compact support contained in $int\, |L|$.
Since $H_c^i(int\, |L|)$
is finitely generated we conclude that $f^*H_c^i(int\, |L|)$ is also 
finitely
generated.  The Claim follows.\\

Now let $r'$ be as in the Claim and let $z_1,\dots , z_l$ be
compactly supported cocycles representing a set of generators of
$\iota_*(H_c^i(B_X(x,r)))$. Since we are assuming $H_c^i(X)=0$, there is
$s=s(r)$ such that $z_1,\dots , z_l$ compactly cobound in $B_X(x,s)$.
Since $X$ is cocompact it follows easily that, given $r$, a number $s$
independent of $x$ exists with the required property.  q.e.d.\\

\noindent {\bf Remark.} By Theorem 1(i) we have $H_c^i(X)=0$, for 
$i>dim\, K$, hence
the number $s$ can be chosen to be independent of $i$.\\

\noindent {\bf Proof of Theorem 1(ii).}   We will give enough 
information for the reader to understand the proof but we will omit 
some details when they are the same as corresponding 
steps in the proof of Theorem B of \cite{O}.  (Theorem B of  \cite{O} 
is the special case of Theorem 1 mentioned in \S 1.)\\ 

By Part (i) it is enough to prove that  $H^*_c(X)$ is 
non-trivial.   Suppose that $H^*_c(X)=0$.
It follows then from Proposition 2.1  that each $H^i_c(X)$ is 
uniformly trivial.  Let $K$ be as in the proof of Part (i).  
Since $X$ is uniformly contractible and $K$ is bounded homotopy 
equivalent to $X$, it follows that $K$ has the following two 
properties:
\begin{enumerate}
\item[{1.}]  $|K|$ is uniformly contractible.

\item[{2.}] $H^*_c(|K|)$ is uniformly trivial (with constants 
independent of the dimension - see Remark after the proof of 
Proposition 2.1).
\end{enumerate}

To obtain a contradiction we now proceed as in the proof of 
Proposition B of \cite{O}.
Embed $|K|$ properly, in some ${\Bbb{R}}^{n}$ by an embedding which 
is affine on each simplex.
Let $T$ be a triangulation of ${\Bbb{R}}^{n}$ such that there is 
a full subcomplex $J$ of $T$
with $| K| = | J|$ and  $J$ is a 
subdivision of $K$. 
Denote the unit metric on $| K|$
by $d_{K}$. As mentioned in \cite{O}, we can extend this 
proper
metric to a proper piecewise
flat metric $d$ on ${\Bbb{R}}^{n}=| T|$.  Indeed, we can assume, 
perhaps after a 
subdivision of $T$ away from $K$, that ${\it mesh}(T) \leq 1$ and 
that every simplex of $|J|$ is convex in $|T|$.\\

Let $M = |N(T',J)|$ where $T'$ is a first derived subdivision of $T$ 
away from $J$, and $N(.,.)$ denotes the simplicial neighborhood.  Then
$M$ is a regular $PL$ neighborhood of $| J| $ 
in  ${\Bbb{R}}^{n}$ \cite{RS}, so
$M$ is an $n$-manifold with boundary $\partial M$ and $M$ is 
bounded homotopy equivalent to $| K| $. \\

Denote by $d_{M}$ the  intrinsic metric on $M$ induced from 
${\Bbb{R}}^{n}$ with metric $d$; i.e. the distance between points of 
$M$ is the greatest lower bound of the $d$-lengths of $PL$ paths in 
$M$ joining them.
Since $|J|\subset M$ we have $d_{M}|_{J}\leq d_{K}$.\\ 

As in \cite{O}, we can choose $M$ close to $| J|$ to get $d_{K}$
close to $d_{M}|_{J}$:\\

\noindent {\bf Claim 1.} We can choose $M$  so that $d_{K}\leq 
(d_{M}|_{J})+1$.\\

The proof is the same as the proof of Claim 1 in the proof of 
Proposition B of \cite{O}.\\

Let $c_t: M\ra |K|$ be the canonical deformation retraction (see 
section 4 of \cite{O}).
As in \cite{O} we also have\\

{\bf (i)}\hspace{.4in} $d_{M}(x,c_t(x))\leq 1$ for every  $t\in 
[0,1]$, and
\vspace{.2in}

{\bf (ii)}\hspace{.2in}  for every $x\in M$ there is $w\in\partial M$ 
with $d_{M}(x,w) \leq 1$.
\vspace{.1in}

Since $|K|$ is bounded homotopy equivalent to $M$ it follows from 
Properties 1 and 2 above that $M$ has the properties:
\begin{enumerate}
\item[{1$'$.}]  $M$ is uniformly contractible.

\item[{2$'$.}] $H^*_c(M)$ is uniformly trivial.
\end{enumerate}

Since $M$ is contractible and 
$H^{i}_{c}(M)=H^{i}_{c}(\, | K| \, )=0$, for all $i$, using the same 
argument as in the
beginning of the proof of Proposition B in \cite{O}, we can arrange 
everything up to this point so  that the pair $(M, \partial M)$ is 
$PL$-homeomorphic to the pair $({\Bbb{R}}^{n}_{+},  {\Bbb{R}}^{n-1})$ 
where
${\Bbb{R}}^{n}_{+}=\{x=(x_{1},...,x_{n})\in {\Bbb{R}}^{n}\,:\,\, 
x_{1}\geq 0\}$.   (See also Remark 2.2 below.)\\

Now, instead of Claims 2 and 3 of \cite{O} (because here we do not 
have a $\Gamma$ action on $K$) we have the following:\\

\noindent {\bf Claim 2.} {\it $\partial M$ is uniformly acyclic, i.e. 
given $b>0$,  there is an $a_{b}$, 
such that, for any $x\in\p M$,  the the homomorphism $H_*(B_{\p M}(x,b) )\ra 
H_*(B_{\p M}(x,a_b))$, induced by the 
inclusion $B_{\p M}(x,b)\hookrightarrow B_{\p M}(x,a_b)$, is zero. 
(Here we consider $\p M$
with the metric $d_M|_{\p M}$.})\\

In other words, every cycle in $\p M$ supported in a ball of radius 
$b$ (with respect to the metric $d_M$) bounds (in $\p M$) in a ball 
of radius $a_b$.\\

\noindent {\bf Proof of Claim 2.} 
Let $z$ be a cycle supported in some ball $B_{\p M}(x,b)$.  By 
Property 1$'$, $M$ is uniformly contractible, 
therefore $z$ bounds in some ball
$B_M(x,b')$, where $b'$ depends only on $b$. In symbols we have $z=\p 
v$, where $v$ is a chain in
$B_M(x,b')$. Since $\p v=z$ and $z$ is a cycle in $\p M$, we have 
that $v$ is a relative cycle in $(M,\p M)$, hence
it is Poincar\'e dual to a cocycle in $M$ that vanishes outside 
$B_M(x,b')$.\\

By Property 2$'$ , $H^*_c(M)$ vanishes uniformly. Hence, there is a 
$b''$ such that the
cocycle dual to $v$ compactly cobounds in $B_M(x,b'')$, again by 
Poincar\'e Duality. Thus there is a chain  $w$ in 
$B_M(x,b'')$ such that $\p w=v+y$, where $y$ is some chain in $\p 
M\cap B_M(x,b'')$.
Then $0=\p\p w=\p (v+y)=\p v+\p y$ and it follows that $z=\p v=\p 
(-y)$. That is, $z$ bounds in
$\p M\cap B_M(x,b'') =  B_{\p M}(x,b'')$.  Take $a_b=b''$. This 
completes the proof of Claim 2.\\

\noindent {\bf Claim 3.} {\it There is a number  $a\geq 0$ such that 
the following holds. Let $P$ be any
subset of $M$ and let $v$ be a relative $i-$cycle of $(M,\p M)$ 
supported in $P$ where $i \leq n-1$ (hence $\p v$ is a cycle in
$\p M\cap P$). Then there is a chain $u$ in $\p M\cap \cN_a(P)$, with 
$\p v=\p u$. (Here
$\cN _a(P)$ denotes the $a$-neighborhood of the set $P$.})\\

The proof of Claim 3 is by induction on the dimension of the 
simplices of the cycle. It uses (i), (ii) and Claim 2 
and is similar to the proof of Claim 4 in the proof of Proposition B 
of \cite{O}.\\

Finally, to obtain a contradiction we use an argument similar to the 
one used at the end
of the proof of Proposition B of \cite{O}:

Let $N=a+2$, where $a$ is the constant from claim 3, and choose a base
point $x_{0}\in \partial M $. Recall that we can assume that $(M,\partial
M)$ is  $PL$-homeomorphic to $({\Bbb{R}}^{n}_{+},{\Bbb{R}}^{n-1})$. Recall
also that $(M,d_{M})$ is a proper metric space.  Let $A$ denote the closed
ball in $M$ with center $x_{0}$ and radius $N$.  Since $A$ is compact,
there is a $PL$ $(n-1)$-ball ${\tilde D}\subset \partial M$, such that
$A\cap\partial M\subset int\, \tilde D$.  Let ${\tilde {\rho}} :S^{n-
2}\rightarrow \partial M $ be a $PL$-embedding such that ${\tilde{\rho}}
(S)=\partial ({\tilde D})$. Since  $x_{0}\in int\, {\tilde D}$ we
have that ${\tilde{\rho}}$ represents a cycle $z$ which does not bound
in $\partial M\setminus\{ x_{0}\}$.  Note  that the cycle $z$ bounds
in $M\setminus A$.  Consequently, by Claim 3, there is a chain  $u$
in $\p M\cap \cN_a(M\setminus A)$ with $\p u=\p v=z$.  But $\p M\cap
\cN_a(M\setminus A)\subset\p M\setminus \{ x_0\}$, which implies that $z$
bounds in $\p M\setminus \{ x_0\}$. This is a contradiction. q.e.d.\\

\noindent {\bf Remark 2.2.} There is a step in the proof where (by 
referring to \cite{O}) we use Stallings' Theorem (\cite{St}) that a
contractible PL manifold of high dimension $d$ which is simply 
connected at infinity is PL homeomorphic to ${\Bbb{R}}^{d}$ .  And to 
apply
that, some tricks are used in \cite{O} to get $\partial M$ to be 
simply connected at infinity.  Stallings' Theorem requires knowledge 
of
engulfing, something invented in order to prove the high-dimensional 
Poincar\'e Conjecture, and certainly more complicated than what ought 
to
be necessary here.  Examination of the above proof, however, shows 
that all we actually need is that the pro-homology at the end of  
$\partial M$ 
in dimension $n-2$ is stably $\Bbb{Z}$, something which follows from 
Poincar\'e Duality.   
The interested reader can supply the details of this; for background, 
see, for example, \cite{G_1}.\\

We denote the Alexander-Spanier cohomology  of the pair $(Z,A)$ by 
$\h^n(Z,A)$.\\

For the proof of Theorem 6, we need a Lemma:\\

{\bf Lemma 2.3.} {\it Let $Z$ be a compact metric space and $A$ a 
closed subset of $Z$. Assume that $\h^n(Z,A)\neq 0$ and that
$\h^{n+1}(Z,B)=0$ for every closed subset $B$ of $Z$. Then there is a 
sequence of open balls $\{ B_k\}_{k\geq k_0}$ in $Z\sm A$, $B_k$ of 
radius $1/k$, such that the homomorphisms
$\h^n(Z,Z\sm B_k )\ra\h^n(Z,A)$ are non-zero. Moreover, we can choose 
the $B_k$'s to satisfy
$d_Z( B_k,A )\geq\delta$, for some $\delta>0$.}\\

{\bf Proof.} Since $\h^*(Z,A)=\h^*_c(Z\sm A)$ (see 6.6.11 of 
\cite{Sp}), elements can be represented
by cocycles with compact support in $Z\sm A$. Hence we can find a 
closed set $A'$, such that $A\sbs int\, A'$ and
the morphism $\h^n(Z,A')\ra\h^n(Z,A)$ is non-zero. Let $k_0$ be such that 
$2/k_0<d_Z(Z\sm A',A)$.\\

Let $U_1,\dots,U_j$ be a finite cover of $\overline{Z\sm A'}$  by 
balls of radius $1/k_0$. Write $U=U_1$ and $V=\bigcup_{i=2}^{\, 
j}U_i$.
Note that $A\sbs Z\sm (U\cup V)\sbs A'$.
We have the following diagram of Mayer-Vietoris sequences:\\

{\scriptsize $$\begin{array}{ccccc}\h^n(Z, Z\sm U)\oplus \h^n(Z, Z\sm 
V)&\ra& \h^n(Z, Z\sm (U\cup V))&\ra&\h^{n+1}(Z, Z\sm (U\cap V))\\ \\
\downarrow\,\,\,\,\,\,\,\,\,\,\,\,\,\,\,\,\,\,\,\,\,\,\,\downarrow&&
\downarrow&&\downarrow 
\\ \\
\h^n(Z,A)\oplus \h^n(Z,A)&\ra& \h^n(Z, A)&\ra&\h^{n+1}(Z, A)
\end{array}$$}
\vspace{.2in}

\noindent where the last group of the first row is zero, by
hypothesis.   Since the non-zero morphism $\h^n(Z,A')\ra\h^n(Z,A)$
factors through $\h^n(Z, Z\sm (U\cup V))$, the middle vertical morphism in
the diagram is non-zero. It follows from the diagram that
(at least) one of the morphisms $\h^n(Z, Z\sm U)\ra\h^n(Z, A)$ and $\h^n(Z,
Z\sm V)\ra\h^n(Z, A)$ is non-zero. If the first one is non-zero take
$B_{k_0}=U$. Otherwise write $U=U_2$ and $V=\bigcup_{i=3}^{\, j}U_i$,
and repeat the process (using the fact that the latter morphism is
non-zero).  Eventually, we will find an $i$ such that the morphism $\h^n(Z,
Z\sm U_i )\ra\h^n(Z, A)$ is non-zero. Take $B_{k_0}=U_i$.  To find
$B_{k_0+1}$ take a finite cover $U_1,\dots,U_j$ of $\bar{ B}_{k_0}$  by
balls of radius $1/k_0+1$ and repeat the process. In this way we obtain
a sequence of balls $B_k$. Note that $d_Z(B_k,A)\geq d_Z(A,\overline{Z\sm
A'})-1/k_0 >0$.  Take $\delta =d_Z(A,\overline{Z\sm A'})-1/k_0$. q.e.d.\\

We recall some facts of dimension theory; see \cite{HW} for details.
Let $Z$ be a compact metric space.  The (Lebesgue covering) {\it
dimension} of $Z$, $dimZ$, is $\leq n$ if every open cover of $Z$ has 
a
refinement whose nerve has dimension at most $n$; one writes $dimZ = 
n$
if, in addition, $dimZ$ is not $\leq n-1$.  If there is no such $n$ 
then
${\it dim} Z$ is infinite.  The $\Bbb{Z}$-{\it cohomological 
dimension}
of $Z$, ${\it dim}_{\Bbb{Z}}Z$, is $n$ if $\h^n(Z,A)\neq 0$ for some
closed subset $A$ of $Z$ while for all $k > n$ and all closed subsets 
$B$
of $Z$, $\h^k(Z,B)= 0$. If there is no such $n$ then   
$dim_{\Bbb{Z}}Z$
is infinite.  Traditionally, here, $\h^n$ refers to {\v C}ech 
cohomology,
but that  is canonically isomorphic to Alexander-Spanier cohomology,
at least when $dimZ$ is finite, as it is in our case - see page 342 of
\cite{Sp}.  There is always the inequality $dim_{\Bbb{Z}}Z \leq dimZ$,
and equality holds when $dimZ$ is finite. \\

{\bf Proof of Theorem 6.}

We write $\bX =X\cup \partial _{\infty} X$.  By Corollary 2, 
$dim\partial _{\infty} X \geq d-1$. We are to prove equality.  By 
Theorem 0 we know that $dim\partial _{\infty} X$ is finite so we may 
use the cohomological definition of dimension.
Suppose that $dim_{\Bbb{Z}} \partial _{\infty} X$ is $n \geq d$. Then 
there
is a closed set $A$ such that $\h^n(\partial _{\infty} X, A)\neq 0$.
Let $B_k$ and $\delta >0$ be as in Lemma 2.3 above (taking 
$Z=\partial _{\infty} X$ with a metric $d_{\partial _{\infty} X}$ 
that induces the cone topology on $\partial _{\infty} X$). Moreover, 
we may assume that the balls $B_k$ converge to some point 
$\gamma\notin A$.\\

Fix a base point $x_0\in X$. For any set $G\sbs \partial _{\infty} 
X$,  the cone $CG$ of $G$ is the union of all
rays emanating from $x_0$ and ending in $G$. Using geodesic 
retraction along rays emanating from $x_0$ we can find a
closed set $D\sbs \bX$ such that:\\

(i) $D\cap \partial _{\infty} X = A$.

(ii) $CA\sbs D$ and $D\sm A$ lies in the 1-neighborhood (in $X$) of 
$CA\sm A$.

(iii) $D$ is a strong deformation retract of $\bX$.\\

By Theorem 1 $H^{n+1}_c(X)$ is trivial, hence uniformly trivial by 
Proposition 2.1, so there is a number $s$ such that every
$(n+1)$-cocycle $y$ with compact support of diameter less than 1 
cobounds a cochain whose compact support lies in the $s$-neighborhood 
of the support of $y$.\\

Since $\gamma\notin A$ we can find a $x_1$ in $[x_0,\gamma]$, such 
that
the ball $B_X(x_1, s+2)$ does not intersect $D$. Write 
$a=d_X(x_0,x_1)$.
Since $B_k\ra \gamma$ there is a $k'$ such that every geodesic ray in
$C\bar{B}_{k'}$ intersects $B_X(x_1, 1)\cap S_X(x_0,a)$ where 
$S_X(x_0,a)$ denotes the sphere of radius $a$ centered at $x_0$ . 
Write $B'=B_{k'}$. By Lemma 2.3,
$\h^n(\partial _{\infty} X,\partial _{\infty} X\sm B' 
)\ra\h^n(\partial
_{\infty} X,A)$ is non-zero.  Let $\{z\}\in\h^n(\partial _{\infty}
X,\partial _{\infty} X\sm B' )$ be such that its image in 
$\h^n(\partial
_{\infty} X,A)$ is non-zero, where $z$ denotes a cocycle with
compact support lying in $B'$.  (Here, as at the beginning of the 
proof of Lemma 2.3, we
are identifying a relative {\v C}ech cohomology group with the 
compactly
supported cohomology of the complement - we will use this convention
again below.)   Using
a geodesic deformation retraction we can find a closed set $E\sbs \bX$
such that:\\

(i) $E\cap \partial _{\infty} X =\partial _{\infty} X\sm B'$.
\vspace{.1in}

(ii) $ C(\partial _{\infty} X\sm B')\cup B_X(x_0, a)\sbs E$.\\

From the exact sequence of the triple $(\bX , \partial _{\infty} 
X\cup D , D)$ and the fact that $\h^*(\bX, D)=0$
we conclude that $\h^n(\partial _{\infty} X\cup D, D )\ra\h^{n+1}(\bX 
,\partial _{\infty} X\cup D)$ is an isomorphism. 
By considering also the exact sequence of the triple $(\bX , \partial 
_{\infty} X\cup E , E)$ and the inclusion
$(\bX , \partial _{\infty} X\cup D , D)\hookrightarrow (\bX , 
\partial _{\infty} X\cup E , E)$ we get the following commutative 
diagram:\\

$$\begin{array}{ccccc}\h^n(\partial _{\infty} X, \partial _{\infty} 
X\sm B')&\cong& \h^n(\partial _{\infty} X\cup E, E)&\ra&\h^{n+1}(\bX, 
\partial _{\infty} X\cup E )\\ \\
\downarrow&&\downarrow&&\downarrow \\ \\
\h^n(\partial _{\infty} X, A)&\cong& \h^n(\partial _{\infty} X\cup D, 
D)&\cong&\h^{n+1}(\bX, \partial _{\infty} X\cup D )
\end{array}$$
\vspace{.2in}

\noindent where the isomorphisms on the left are given by excision
- see 6.6.5 of \cite{Sp}. Let $\{z'\}$ be the image of $\{z\}$ in
$\h^{n+1}(\bX, \partial _{\infty} X\cup E )$. Then the image of $\{z'\}$
in $\h^{n+1}(\bX, \partial _{\infty} X\cup D )\cong \h_c^{n+1}(X,D\sm A)$
is non-zero.  We regard $z'$ as a cocycle compactly supported outside
$\partial _{\infty} X\cup E$ and we obtain a contradiction by showing
that $z'$ compactly cobounds outside $D\sm A$ as follows: using a proper
radial contraction along rays emanating from $x_0$ we find a cochain $c$
representing $z'$ with compact support lying in $B_X(x_1, 1)$. This $z'$
cobounds in $B_X(x_1, s+1)$ - a contradiction since $B_X(x_1, s+1)$
is disjoint from $D$. q.e.d.\\

Ross Geoghegan, Binghamton University (SUNY), Binghamton  N.Y., 13902,
U.S.A.\\

Pedro Ontaneda, UFPE, Recife, PE 50670-901, Brazil


\begin{thebibliography}{99}


\bibitem{BG} R. Bieri and R.Geoghegan,
{\em Connectivity properties of group actions on
non-positively curved spaces}, Mem. Amer. Math. Soc. vol. 161, no. 
{\bf 765} (2003).

\bibitem{Be} M. Bestvina,  {\em Local homology properties of 
boundaries of groups.} Michigan Math. J. 43 (1996), no. 1, 123--139.

\bibitem{BM} M. Bestvina and G. Mess,  {\em The boundary of 
negatively curved groups} J. Amer. Math. Soc. 4 (1991), no. 3, 
469--481.

\bibitem{BH} M. Bridson and A. Haeflinger, {\em Metric spaces of
non-positive curvature}, Grundlehren der Mathematischen Wissenschaften
{\bf 319}. Springer-Verlag, Berlin, (1999).

\bibitem{Br} K.S. Brown, {\em Cohomology of groups} Graduate Texts in 
Mathematics, 87. Springer-Verlag, New York-Berlin, (1982)

\bibitem{CK} C. Croke and B. Kleiner {\em Spaces with nonpositive 
curvature and their ideal boundaries} Topology 39 (2000), no. 3, 
549--556

\bibitem{G_1} R. Geoghegan,  {\em Topological methods in group 
theory},  monograph in
preparation.

\bibitem{G_2} R. Geoghegan, {\em  The shape of a group--connections 
between shape theory
and the homology of groups},  Geometric and algebraic topology, 
271--280,
Banach Center Publ., {\bf 18}, PWN, Warsaw  (1986).
 
\bibitem{GM} R. Geoghegan and M. Mihalik, {\em Free abelian 
cohomology of groups
and ends of universal covers}, J. Pure Appl. Algebra 36 (1985)
123-137.

\bibitem{Ger} S. Gersten,  {\em Quasi-isometry invariance of 
cohomological
dimension} C. R. Acad. Sci. Paris  Sci. Paris Sci. Paris Sér. I Math. 
316 (1993),
no. 5, 411--416.

\bibitem{GH} M. J. Greenberg and J. R. Harper, {\em Algebraic 
topology: a first
course},  Math. Lecture Notes Series {\bf 58} Benjamin/Cummings 
Publishing Co., Inc. (1981).

\bibitem{HW} W. Hurewicz and H. Wallman, {\em Dimension theory}
Princeton Mathematical Series, vol. 4. Princeton University Press,
Princeton, N. J., (1941).

\bibitem{K} B. Kleiner,  {\em The local structure of length spaces 
with
curvature bounded above} Math. Z. 231 (1999), no. 3, 409--456.

\bibitem{O} P. Ontaneda, {\em Cocompact CAT(0) spaces are almost
geodesically complete}.  Topology (in press).

\bibitem{RS} C. Rourke and D. Sanderson, {\em Introduction to 
piecewise
linear topology}. Ergebnisse der Mathematik und ihrer Grenzgebiete,
Band 69. Springer-Verlag, New York-Heidelberg, (1972).
                                                             
\bibitem{Sp} E. H. Spanier, {\em Algebraic topology}  McGraw-Hill 
Book Co., New
York-Toronto-London (1966). 

\bibitem{St} J. Stallings, {\em The piecewise-linear structure of 
euclidean
space}, Proc. Cambridge Philos. Soc. 58 (1962) 481-488.

\bibitem{Sw} E. L. Swenson, {\em A cut point theorem for CAT(0) 
groups}, Jour. Diff. Geom. 53 (1999) 327-358.

\bibitem{W} A. Weil, Sur les th\'eor\`emes de de Rham. Comment. Math. 
Helv. 26,
(1952). 119--145.


\end{thebibliography}
\end{document}